\documentclass[11pt]{article}
\usepackage{amsfonts}
\usepackage{amssymb}
\usepackage{amscd}
\usepackage{amsmath}
\usepackage{epsfig}
\usepackage{graphicx, subfigure}
\usepackage{latexsym}

\newlength{\dinwidth}
\newlength{\dinmargin}
\newtheorem{theorem}{Theorem}

\newtheorem{corollary}{Corollary}
\newtheorem{remark}{Remark}

\begin{document}

\title{Note on algebro-geometric solutions to triangular  Schlesinger systems}
\author{Vladimir Dragovi\'c$^1$ and Vasilisa Shramchenko$^2$}

\date{}

\maketitle

\footnotetext[1]{Department of Mathematical Sciences, University
of Texas at Dallas, 800 West Campbell Road, Richardson TX 75080,
USA. Mathematical Institute SANU, Kneza Mihaila 36, 11000
Belgrade, Serbia.  E-mail: {\tt
Vladimir.Dragovic@utdallas.edu}; corresponding author.}

\footnotetext[2]{Department of mathematics, University of
Sherbrooke, 2500, boul. de l'Universit\'e,  J1K 2R1 Sherbrooke, Quebec, Canada. E-mail: {\tt Vasilisa.Shramchenko@Usherbrooke.ca}
%--the corresponding author
}

\begin{abstract}

We construct algebro-geometric upper triangular solutions of rank two Schlesinger systems. Using these solutions we derive two families  of  solutions to the sixth Painlev\'e equation with parameters $({1}/{8}, -{1}/{8}, {1}/{8}, {3}/{8})$  expressed in simple forms using periods of differentials on elliptic curves. Similarly  for every integer $n$ different from $0$ and $-1$ we obtain one family of  solutions to the sixth Painlev\'e equation with parameters $(\frac{9n^2+12n+4}{8}, -\frac{n^2}{8}, \frac{n^2}{8}, \frac{4-n^2}{8})$.

\end{abstract}

\vskip 1cm
MSC primary 34M55, 34M56; secondary 14H70.
\vskip 1cm

\section{Introduction. Schlesinger system}

Consider the Fuchsian matrix linear system for a $2\times 2$ matrix function $\Phi(u)$ defined on the Riemann sphere
\begin{equation*}
\label{linsys_intro}
\frac{d\Phi}{du} = A(u) \Phi, \qquad u \in {\mathbb CP}^1
\end{equation*}
where the matrix $A\in sl(2,{\mathbb C})$ has $2g+2$ simple poles at the branch points of a hyperelliptic curve given by the equation $v^2=(u-u_1)\cdots(u-u_{2g+1})$, that is at  $u_1, u_2, \dots,
u_{2g+1}\in{\mathbb CP}^1$ and at the point at infinity:
\begin{equation}
\label{Ah}
A(u) = \sum_{j=1}^{2g+1} \frac{A^{(j)}}{u-u_j}.
\end{equation}
The equation of isomonodromic deformation of the linear system form the
Schlesinger system of partial differential equations for the residue-matrices $A^{(j)}\in
sl(2,{\mathbb C})$ with respect to positions $u_i$ of poles of $A$
as independent variables:

\begin{equation}
\label{hSchlesinger_intro}
\frac{\partial A^{(j)}}{\partial u_k} =
\frac{[A^{(k)}, A^{(j)}]}{u_k-u_j}; \qquad\qquad  \frac{\partial
A^{(k)}}{\partial u_k} = -\sum_{j\neq k}\frac{[A^{(k)},
A^{(j)}]}{u_k-u_j},
\end{equation}
where  $ A^{(\infty)}:=-A^{(1)} -\dots-A^{(2g+1)}  = const.$

In the simplest nontrivial case $g=1$, the associated curve is elliptic and is usually given in the Legendre form $v^2=u(u-1)(u-x)$.
The Schlesinger system then reads
\begin{equation}
\label{Schlesinger}
\frac{dA^{(1)}}{dx} = \frac{[A^{(3)}, A^{(1)}]}{x}; \qquad
\frac{dA^{(2)}}{dx} = \frac{[A^{(3)}, A^{(2)}]}{x-1}; \qquad
A^{(\infty)}:= -A^{(1)} -A^{(2)}-A^{(3)} = const.
\end{equation}
This Schlesinger system is well known to be closely related to (``equivalent to'', see \cite{JimboMiwa}, \cite{JMU}) the Painlev\'e VI equation.
The Painlev\'e VI equation is the second order ordinary differential equation
\begin{equation}
\label{Painleve}
\frac{d^2 y}{dx^2} = \frac{1}{2} \left(  \frac{1}{y} + \frac{1}{y-1}+ \frac{1}{y-x}  \right) \left( \frac{dy}{dx} \right)^2 - \left( \frac{1}{x} + \frac{1}{x-1} + \frac{1}{y-x} \right) \frac{dy}{dx}
\end{equation}
\begin{equation*}
 + \frac{y(y-1)(y-x)}{x^2(x-1)^2}\left( \hat \alpha +\hat \beta\frac{x}{y^2} + \hat \gamma\frac{x-1}{(y-1)^2} + \hat \delta\frac{x(x-1)}{(y-x)^2}  \right)
\end{equation*}
with parameters $\hat \alpha, \hat \beta, \hat \gamma, \hat \delta \in \mathbb C$.

Using the freedom of  global conjugation of the residue matrices by a constant invertible matrix, one may assume  $A^{(\infty)}$ to be diagonal.  Then the entry $A_{12}(u)$ of the matrix $A(u)$  is of the form
\begin{equation*}
\label{A12}
A_{12}(u) = \kappa\frac{(u-y)}{u(u-1)(u-x)}
\end{equation*}
where $\kappa$ is a function of $x$. The position of the only zero of $A_{12}(u)$ as a function of the position $x$ of the pole is the function $y(x)$ which satisfies
the Painlev\'e VI equation (\ref{Painleve}) with parameters
\begin{equation}\label{constants1}
\hat\alpha = \frac{(2\alpha_\infty-1)^2}{2}, \qquad \hat\beta=-2\alpha_1^2, \qquad \hat\gamma=2\alpha_2^2, \qquad \hat\delta = \frac{1}{2} - 2\alpha_3^2.
\end{equation}
Here, $\alpha_i$ and $-\alpha_i$ are the eigenvalues of the residue-matrix $A^{(i)}$; these quantities are integrals of motion of system (\ref{Schlesinger}).

Although the above Schlesinger system for $g=1$ implies
that the corresponding function $y(x)$ solves an appropriate Painlev\'e VI equation, the correspondence between the Schlesinger system and the Painlev\'e equation is not one to one. This is because the relation (\ref{constants1}) is nonlinear. Moreover, the Schlesinger system admits a simple, but nontrivial reduction to a triangular case, see \cite{book} and \cite{GontsovLeksin}.

The main purpose of this note is to construct algebro-geometric upper triangular solutions of the above Schlesinger systems. There are two well known approaches to algebro-geometric solutions to the Schlesinger systems, both proposed almost twenty years ago in \cite{Deift} and \cite{KiKo}. Recently, the authors of this note presented yet another
approach, see \cite{DS1, DS}. However, in each of the three papers, the triangular reduction was left out of the scope. On the other hand, in each of the three papers, constructed solutions to
the Schlesinger systems had eigenvalues $\alpha_i = \pm\frac{1}{4}$.  In these papers in the basic $g=1$ case the main role was attributed to the Painlev\'e VI equation (\ref{Painleve})  with the parameters
\begin{equation}
\label{constants}
\hat \alpha=\frac{1}{8}, \qquad \hat \beta=-\frac{1}{8}, \qquad \hat\gamma=\frac{1}{8}, \qquad \hat\delta=\frac{3}{8}.
\end{equation}
In the present  note, we consider  upper triangular Schlesinger systems, for all of which  the eigenvalues of the residue matrices $A^{(i)}$ are  $\pm n/4$ with nonzero integer $n$. In the case $g=1$ and $n=\pm1$, two of our Schlesinger systems give rise to  the Painlev\'e VI equation  (\ref{Painleve}) with the parameters (\ref{constants}).
Thus, we find two families of solutions of PVI$(1/8, -1/8, 1/8, 3/8)$ in new algebro-geometric forms. These families are included in the general solution, see \cite{Hitchin, KiKo, Okamoto} and Section \ref{sect_others}. In the case $g=1$ and $n$ different from $0$ and $-1$, the Schlesinger systems considered correspond to the Painlev\'e VI equation with the parameters
\begin{equation}
\label{intro_npara}
\left(\frac{9n^2+12n+4}{8}, -\frac{n^2}{8}, \frac{n^2}{8}, \frac{4-n^2}{8}\right)\;.
\end{equation}
Note that upper triangular Schlesinger systems were also studied in \cite{book, GontsovLeksin, Leksin}. In \cite{KapovichMillson, Kohno} solutions similar to ours were obtained in a different context and for a different system of equations, see Section~\ref{sect_triangular}.

\section{Upper triangular solutions of rank two Schlesinger system and Euler-Poisson-Darboux equations}
\label{sect_triangular}

We are looking for solutions to the Schlesinger system in the form:
\begin{equation*}
A^{(i)} = \left(\begin{array}{cc}\alpha_i & a_i \\0 & -\alpha_i\end{array}\right), \qquad i=1,\dots, 2g+1
\end{equation*}
where $\alpha_i$ are constants, $a_i$ are functions of $u_1, \dots, u_{2g+1}$  and
\begin{equation*}
A^{(\infty)} = \left(\begin{array}{cc}\alpha_\infty & 0 \\0 & -\alpha_\infty\end{array}\right).
\end{equation*}

The Schlesinger system becomes a system of equations for the functions $a_i$:
\begin{equation}
\label{Schlesinger_ai}
\frac{\partial a_i}{\partial u_j} = \frac{2(\alpha_j a_i - \alpha_i a_j)}{u_j-u_i}\;\;\;\mbox{ for } \;\; i\neq j, \;\;\mbox{ and } \qquad \sum_{i=1}^{2g+1} a_i=0.
\end{equation}

Note that the system of PDEs in \eqref{Schlesinger_ai} is potential, that is there exists a function $f=f(u_1,\dots, u_{2g+1})$ such that $a_i=\partial_{u_i}f$. Therefore this system (without the condition $\sum_{i=1}^{2g+1} a_i=0$) is the Euler-Poisson-Darboux system
\begin{equation}
\label{EPD}
\frac{\partial^2 f}{\partial u_j \partial u_i} = \frac{1}{u_i-u_j} \left( \beta_j\frac{\partial f}{\partial u_i}  - \beta_i\frac{\partial f}{\partial u_j} \right)
\;\;\;\mbox{ for } \;\; i\neq j,
\end{equation}
where we denote $\beta_i=-2\alpha_i.$ For more on the Euler-Poisson-Darboux system see for example \cite{KodamaKonopelchenko, Konopelchenko} The Euler-Poisson-Darboux equations are also linked to confocal coordinates in $\mathbb R^n$, see for example \cite{BobenkoSuris}.

A similar system to \eqref{Schlesinger_ai} appeared in \cite{KapovichMillson, Kohno} in a different context. In \cite{KapovichMillson} the system was called the hypergeometric equation. The form of that equation is similar to \eqref{Schlesinger_ai} but the second condition $\sum_{i=1}^{2g+1} a_i=0$ is different.

In this section we consider three particular choices of sets of eigenvalues $\{\alpha_i\}$ and construct families of solutions in these cases in terms of periods of meromorphic differentials on associated elliptic curves. In the case of curves of genus one, the eigenvalues considered in Case 1 with $n=-1$ and Case 2 correspond to the Painlev\'e VI equation with parameters $({1}/{8}, -{1}/{8}, {1}/{8}, {3}/{8})$. The eigenvalues of Case 1 with arbitrary integer $n\neq 0$ correspond to  the Painlev\'e VI equation with parameters \eqref{intro_npara}.

\subsection{Case 1}
\label{sect_Schlesinger1}

Let us  put $\alpha_i=n/4$ with $n\in\mathbb Z\setminus\{0\}$. Then we obtain $\alpha_\infty  = -n(2g+1)/4$ and (\ref{Schlesinger_ai}) becomes
\begin{equation}
\label{aisystem}
\frac{\partial a_i}{\partial u_j} = \frac{n}{2}\frac{a_i -  a_j}{(u_j-u_i)} \;\;\;\mbox{ for } \;\; i\neq j, \;\;\mbox{ and } \qquad \sum_{i=1}^{2g+1} a_i=0.
\end{equation}

The Schlesinger system (\ref{hSchlesinger_intro}) is naturally associated with a family of hyperelliptic curves defined by the equation
\begin{equation}
\label{hyper}
v^2=(u-u_1)\cdots(u-u_{2g+1}),
\end{equation}
where the varying branch points $u_1, \dots, u_{2g+1}$ are given by the positions of the poles in (\ref{Ah}).

Consider the differential
\begin{equation*}
\phi=v^n du\,.
\end{equation*}
This is a meromorphic differential on the curve (\ref{hyper}) (the differential is holomorphic if $n=-1$).  Its period $\oint_\gamma \phi$  over a closed cycle $\gamma$  is a function of the branch points of the curve. We build a solution to the Schlesinger system in terms of  derivatives of this period:
\begin{equation*}
\frac{\partial \oint_\gamma \phi}{\partial u_i}= -\frac{n}{2} \oint_\gamma\frac{v^n du}{(u-u_i)}.
\end{equation*}
\begin{theorem}
\label{thm}
The functions
\begin{equation*}
a_i= \oint_\gamma\frac{v^n du}{(u-u_i)}
\end{equation*}
with $n\in\mathbb Z\setminus\{0\}$ satisfy system (\ref{aisystem}).
\end{theorem}

{\it Proof.}
Let us compute the derivative of $a_i$ with respect to $u_j$ for $j\neq i:$
\begin{equation*}
\frac{\partial a_i}{\partial u_j} = \frac{-n}{2}\oint_\gamma\frac{v^n du}{(u-u_i)(u-u_j)} = \frac{-n}{2} \frac{1}{u_i-u_j} \oint_\gamma{v^n du}\left( \frac{1}{u-u_i} - \frac{1}{u-u_j}\right) = \frac{-n}{2}\frac{1}{u_i-u_j} (a_i-a_j).
\end{equation*}
Proof of the second relation in (\ref{aisystem}) is also a straightforward computation starting from the definition of $a_i$ given in the theorem:
\begin{equation*}
\sum_{i=1}^{2g+1}a_i = \oint_\gamma\sum_{i=1}^{2g+1}\frac{v^n du}{(u-u_i)} = \frac{2}{n}\oint_\gamma d\left(v^n\right)
\end{equation*}
which is zero as an integral of an exact differential over a closed cycle.
$\Box$

Since the first homology of our hyperelliptic curve (\ref{hyper}) is spanned by $2g$ cycles $\{\gamma_k\}_{k=1}^{2g}$, we have found a $2g-$parameter family of solutions to the Schlesinger system.

\begin{corollary}
\label{cor_Schlesinger}
For arbitrary parameters $c_1,\dots, c_{2g}$ the following triangular matrices solve (\ref{hSchlesinger_intro})
\begin{equation*}
A^{(i)} = \left(\begin{array}{cc}n/4 & a_i \\0 & -n/4\end{array}\right), \qquad i=1,\dots, 2g+1
\end{equation*}
where $n\in\mathbb Z\setminus\{0\}$ and
\begin{equation*}
a_i=\sum_{k=1}^{2g} c_k \oint_{\gamma_k} \frac{v^n du}{(u-u_i)}.
\end{equation*}
\end{corollary}

\begin{remark}
The corresponding solution of the Euler-Poisson-Darboux equation \eqref{EPD} with $\beta_i=-2\alpha_i=-n/2$ is given by
\begin{equation}
\label{potential}
f=-\frac{2}{n}\oint_\gamma v^n du\;,
\end{equation}
where $\gamma=\sum_{k=1}^{2g} c_k\gamma_k$ and $n\in\mathbb Z\setminus\{0\}.$
\end{remark}

\begin{remark}
We note that the integrands $v^n$ in \eqref{potential} solve the Euler-Poisson-Darboux equation \eqref{EPD} with $\beta_i=-n/2$. This solution with $n =1$, namely
\begin{equation*}
f= \prod_{i=1}^{2g+1}(u-u_i)^{1/2},
\end{equation*}
is a so-called separable solution of the Euler-Poisson-Darboux equation. Such solutions are directly related to the Jacobi coordinates associated with confocal quadrics in $\mathbb R^{2g+1}$, see \cite{BobenkoSuris}.
\end{remark}

\subsection{Case 2}
\label{sect_Schlesinger2}

Consider the family \eqref{hyper} of hyperelliptic curves with the first branch point fixed at zero: $u_1=0.$ This can always be achieved by a conformal transformation in the $u$-sphere. In other words, we consider the curves of the form
\begin{equation}
\label{hyper_0}
v^2=u\prod_{i=2}^{2g+1}(u-u_i).
\end{equation}
The associated triangular reduction of the Schlesinger system has the form \eqref{Schlesinger_ai} where the equations involving derivative with respect to $u_1$ are missing.
\begin{theorem}
Let $\{\gamma_k\}_{k=1}^{2g}$ be a basis in homology of the hyperelliptic curve \eqref{hyper_0} and
let $\gamma = \sum_{k=1}^{2g}c_k\gamma_k$. The following functions
\begin{equation*}
a_1=-\oint_\gamma \frac{du}{v}, \qquad a_i= \oint_\gamma \frac{du}{v} + u_i\oint_\gamma \frac{du}{(u-u_i)v}  \;\;\;\mbox{for}\;\;\;\;  2\leq i\leq 2g+1
\end{equation*}
satisfy the Schlesinger system \eqref{Schlesinger_ai} assuming that $u_1=0$ and
with the eigenvalues of the residue matrices given by
\begin{equation*}
\alpha_1=1/4, \qquad \alpha_i=-1/4 \;\;\;\mbox{for}\;\;\;\;  2\leq i\leq 2g+1\;\;\;\mbox{ and }\;\;\;\; \alpha_\infty=(2g-1)/4.
\end{equation*}
\end{theorem}
{\it Proof.}
The Schlesinger system with given eigenvalues reduces to the following equations for the functions $a_i$
\begin{equation}
\label{aisystem_0}
\frac{\partial a_1}{\partial u_j} = -\frac{a_1 + a_j}{2u_j}, \qquad
\frac{\partial a_i}{\partial u_j} = \frac{a_j -  a_i}{2(u_j-u_i)} \;\;\mbox{ for } \;\; i\neq j, \; 2\leq i,j\leq 2g+1\;\;\mbox{ and } \;\;\; \sum_{i=1}^{2g+1} a_i=0.
\end{equation}
The two differential equations are verified by a straightforward computation similarly to Theorem \ref{thm}. To obtain the sum of $a_i$ consider the following integral of an exact form (here we use $u_1$ as indeterminate which will later be set equal to zero):
\begin{equation*}
0=\oint_\gamma d\left(\frac{u}{v}\right) =\oint_\gamma \frac{du}{v} - \frac{1}{2}\oint_\gamma\sum_{i=1}^{2g+1}\frac{u}{u-u_i}\frac{du}{v} =
\oint_\gamma \frac{du}{v} - \frac{1}{2}\oint_\gamma\sum_{i=1}^{2g+1}\left(1+\frac{u_i}{u-u_i}\right)\frac{du}{v} .
\end{equation*}
This implies the relation
\begin{equation}
\label{zero_rel}
(2g-1)\oint_\gamma \frac{du}{v} +\sum_{i=1}^{2g+1} u_i \oint_\gamma \frac{du}{(u-u_i)v}=0.
\end{equation}
Putting $u_1=0$ we obtain the desired $\sum_{i=1}^{2g+1} a_i=0$.
$\Box$

The result of this subsection is summarized in the following corollary.
\begin{corollary}
Let $\{\gamma_k\}_{k=1}^{2g}$ be a basis in homology of the hyperelliptic curve \eqref{hyper_0} and $\gamma = \sum_{k=1}^{2g}c_k\gamma_k$ with arbitrary parameters
$c_1,\dots, c_{2g}$ be a closed cycle.
The following triangular matrices solve (\ref{hSchlesinger_intro})
\begin{equation*}
A^{(1)} = \left(\begin{array}{cc}1/4 & a_1 \\0 & -1/4\end{array}\right), \qquad
A^{(i)} = \left(\begin{array}{cc}-1/4 & a_i \\0 & 1/4\end{array}\right), \qquad i=2,\dots, 2g+1
\end{equation*}
with
\begin{equation*}
a_1=-\oint_\gamma \frac{du}{v}, \qquad a_i= \oint_\gamma \frac{du}{v} + u_i\oint_\gamma \frac{du}{(u-u_i)v}, \qquad  i=2,\dots, 2g+1.
\end{equation*}
\end{corollary}

\begin{remark}
The corresponding solution of the Euler-Poisson-Darboux equation \eqref{EPD} with $\beta_i=-2\alpha_i$ is given by
\begin{equation*}
f=2 \oint_\gamma \frac{u du}{v}\;,
\end{equation*}
where $\gamma=\sum_{k=1}^{2g} c_k\gamma_k.$
\end{remark}

\section{Tau-function}

The isomonodromic tau-function of the Schlesinger system is defined  by
\begin{equation*}
\frac{\partial {\rm ln}\tau}{\partial u_j} = \frac{1}{2}\underset{u=u_j}{\rm res}{\rm tr} A^2(u)du
\end{equation*}
where the matrix $A(u)$ is given by (\ref{Ah}).

In the case of Schlesinger system from Section \ref{sect_Schlesinger1}, for the matrices $A_i$ from Corollary \ref{cor_Schlesinger}, this definition becomes
\begin{equation*}
\frac{\partial {\rm ln}\tau}{\partial u_j} = \frac{n^2}{8}\sum_{i=1, i\neq j}^{2g+1} \frac{1}{u_j-u_i}.
\end{equation*}
Therefore, the tau-function in this case, up to an arbitrary factor, is the following rational power of the Vandermonde determinant (recall that $n\in\mathbb Z\setminus\{0\}$):
\begin{equation*}
\tau=\left(\prod_{i,k=1, k<i}^{2g+1} (u_i-u_k)\right)^{\frac{n^2}{8}}.
\end{equation*}
In the case of Schlesinger system from Section \ref{sect_Schlesinger2}, the definition of the tau-function gives
\begin{equation*}
\frac{\partial {\rm ln}\tau}{\partial u_j} = -\frac{1}{8u_j}  + \frac{1}{8}\sum_{i=2, i\neq j}^{2g+1} \frac{1}{u_j-u_i} \;,\qquad j\geq 2
\end{equation*}
and the tau-function, up to an arbitrary factor, is given by
\begin{equation*}
\tau=\left(\frac{{\underset{i,k=2, k<i}{\overset{2g+1}{\prod}} }(u_i-u_k)}{\prod_{k=2}^{2g+1}u_k}\right)^{\frac{1}{8}}.
\end{equation*}

In particular we see that in both cases the tau-functions do not vanish if all branch points of the curve are distinct and thus the Malgrange divisor is empty.

\section{Examples in genus one: Painlev\'e-VI}

In the elliptic case one usually applies a conformal transformation in the $u$-sphere to map two branch points to $0$ and $1$ and to keep the branch point at infinity. Thus we have $u_1=0$, $u_2=1$ and $u_3=x$.
The case of the elliptic family $v^2=u(u-1)(u-x),$ with only one variable branch point, is included in both cases of the upper triangular Schlesinger system considered above.
As is known, see for example \cite{book}, in the case of this family, the Schlesinger system corresponds to the Painlev\'e VI equation with the parameters computed from the eigenvalues $\alpha_1, \alpha_2, \alpha_3, \alpha_\infty$ of the matrices solving the Schlesinger system as in (\ref{constants1}).

Two of the Schlesinger systems we consider correspond to the Painlev\'e-VI equation with parameters $(\frac{1}{8}, -\frac{1}{8}, \frac{1}{8}, \frac{3}{8})$. These are the systems with eigenvalues $\alpha_1=\alpha_2=\alpha_3=-1/4, \; \alpha_\infty=3/4$  (Section \ref{sect_Schlesinger1}, $n=-1$) and  $\alpha_1=1/4, \;\alpha_2=\alpha_3=-1/4, \; \alpha_\infty=1/4$ (Section \ref{sect_Schlesinger2}). For the system of Section \ref{sect_Schlesinger1} with $n$ different from $0$ and $-1$, that is $\alpha_1=\alpha_2=\alpha_3=n/4$ and $\alpha_\infty=-3n/4$, we obtain the the Painlev\'e-VI equation with parameters \eqref{intro_npara}.

\subsection{Case 1}
\label{sect_Painleve1}

In the case of the elliptic curve $v^2=u(u-1)(u-x),$  our solution to the Schlesinger system with eigenvalues $\alpha_1=\alpha_2=\alpha_3=-1/4, \; \alpha_\infty=3/4$  from Section \ref{sect_Schlesinger1} ($n=-1$) reads
\begin{equation*}
a_1(x)=\oint_\gamma \frac{du}{uv}, \qquad a_2(x)=\oint_\gamma \frac{du}{(u-1)v}, \qquad a_3(x)=\oint_\gamma \frac{du}{(u-x)v}
\end{equation*}
with $\gamma=c_1{\mathcal A}+c_2{\mathcal B}$ where ${\mathcal A}$ and ${\mathcal B}$ are two generators of the homology group of the torus.

The only zero of $A_{12}(u)$, the top right entry of the matrix \eqref{Ah} gives a solution of the associated Painlev\'e equation. Using this fact we can write the solution
 in terms of periods on elliptic curves:
\begin{equation*}
A_{12}(u) = \frac{(u-1)(u-x)\oint_\gamma \frac{du}{uv}+u(u-x)\oint_\gamma \frac{du}{(u-1)v}+u(u-1)\oint_\gamma \frac{du}{(u-x)v}}{u(u-1)(u-x)}.
\end{equation*}
Using the relation $a_1(u)+a_2(u)+a_3(u)=0$, we see that the numerator is linear in $u$:
\begin{equation*}
A_{12}(u) = \frac{x\oint_\gamma \frac{du}{uv} + \left(x\oint_\gamma \frac{du}{(u-x)v}+\oint_\gamma \frac{du}{(u-1)v}\right) u}{u(u-1)(u-x)}.
\end{equation*}
Thus, we obtain a family of solutions of PVI$(\frac{1}{8}, -\frac{1}{8}, \frac{1}{8}, \frac{3}{8})$ in the form:
\begin{equation*}
y(x) = -\frac{x\oint_\gamma \frac{du}{uv}}{ x\oint_\gamma \frac{du}{(u-x)v}+\oint_\gamma \frac{du}{(u-1)v}}
\end{equation*}
where as before $\gamma=c_1{\mathcal A}+c_2{\mathcal B}$. Note that this family of solutions is parametrized by one constant $c_1/c_2$.

This solution can be easily differentiated using equations \eqref{aisystem} of the Schlesinger system. These equations in our example become
\begin{equation*}
\frac{da_1}{dx} = \frac{1}{2}\frac{a_3-a_1}{x}, \qquad \frac{da_2}{dx} = \frac{1}{2}\frac{a_3-a_2}{x-1}, \;\;\mbox{ and }\;\; a_1+a_2+a_3=0.
\end{equation*}
Using them we can rewrite the solution, for example, as follows:
\begin{equation*}
y(x) = \frac{x\oint_\gamma \frac{du}{uv}}{ x\oint_\gamma \frac{du}{uv}+(x-1)\oint_\gamma \frac{du}{(u-1)v}} = \frac{xa_1}{xa_1+(x-1)a_2}
\end{equation*}
and find its first two derivatives:
\begin{equation*}
y'(x) = \frac{xa_1^2 - (x-1)a_2^2}{2(xa_1+(x-1)a_2)^2} \;\;\;\;\;\mbox{ and }\;\;\;\;\;
y''(x) = \frac{-a_1a_2a_3}{2(xa_1+(x-1)a_2)^3}\;.
\end{equation*}
With these expressions, the validity of Painlev\'e-VI with parameters $(\frac{1}{8}, -\frac{1}{8}, \frac{1}{8}, \frac{3}{8})$ is verified by  a short straightforward calculation. Moreover, the expression for $y'$ can be used to compute Okamoto transformations of this solution and thus to express solutions to other Painlev\'e VI equations in terms of elliptic periods.

\begin{remark}
Analogously, for the Schlesinger system from Section \ref{sect_Schlesinger1} with integer $n\notin\{ 0,-1\}$ associated with the curve $v^2=u(u-1)(u-x)$, we obtain a family of solutions of the Painlev\'e-VI equation with parameters $(\frac{9n^2+12n+4}{8}, -\frac{n^2}{8}, \frac{n^2}{8}, \frac{4-n^2}{8})$:
\begin{equation*}
y(x) = -\frac{x\oint_\gamma \frac{v^n du}{u}}{ x\oint_\gamma \frac{v^n du}{u-x}+\oint_\gamma \frac{v^n du}{u-1}}=\frac{x\oint_\gamma \frac{v^ndu}{u}}{ x\oint_\gamma \frac{v^ndu}{u}+(x-1)\oint_\gamma \frac{v^ndu}{u-1}}
\end{equation*}
indexed by the value of  $c_1/c_2$ such that $\gamma=c_1{\mathcal A}+c_2{\mathcal B}$.
We also obtain analogously
\begin{eqnarray*}
&&y'(x) = \frac{  n(x-1)a_2^2-nxa_1^2 - 2(n+1) a_1a_2}{2(xa_1+(x-1)a_2)^2}
\end{eqnarray*}
with
\begin{equation*}
a_1(x)=\oint_\gamma \frac{v^ndu}{u}, \qquad a_2(x)=\oint_\gamma \frac{v^ndu}{u-1}, \qquad a_3(x)=\oint_\gamma \frac{v^ndu}{u-x}\,.
\end{equation*}

\end{remark}

\subsection{Case 2: the Picard-Fuchs equation}

For the same elliptic curve $v^2=u(u-1)(u-x),$  the solution to the Schlesinger system with eigenvalues $\alpha_1=1/4, \;\alpha_2=\alpha_3=-1/4, \; \alpha_\infty=1/4$ from Section \ref{sect_Schlesinger2} reads
\begin{equation*}
a_1(x)=-\oint_\gamma \frac{du}{v}, \qquad a_2(x)=\oint_\gamma \frac{du}{v}+ \oint_\gamma \frac{du}{(u-1)v}  , \qquad a_3(x)=\oint_\gamma \frac{du}{v}+x\oint_\gamma \frac{du}{(u-x)v}
\end{equation*}
with $\gamma=c_1{\mathcal A}+c_2{\mathcal B}$ where ${\mathcal A}$ and ${\mathcal B}$ are two generators of the homology group of the torus.

Using relation \eqref{zero_rel} which in genus one becomes
\begin{equation*}
\oint_\gamma \frac{du}{v} + \oint_\gamma \frac{du}{(u-1)v} +x\oint_\gamma \frac{du}{(u-x)v}=0,
\end{equation*}
we can rewrite this solution in a slightly simpler form
\begin{equation}
\label{Picard-Schlesinger}
a_1(x)=-\oint_\gamma \frac{du}{v}, \qquad a_2(x)=-x\oint_\gamma \frac{du}{(u-x)v}, \qquad a_3(x)=-\oint_\gamma \frac{du}{(u-1)v}.
\end{equation}

Again, a family of solutions of the associated Painlev\'e equation is obtained as the only zero of $A_{12}(u)$ in the same way as in Section \ref{sect_Painleve1}:
\begin{equation*}
y(x)=-\frac{xa_1}{xa_3+a_2} = -\frac{\oint_\gamma \frac{du}{v}}{ \oint_\gamma \frac{du}{(u-1)v}+ \oint_\gamma \frac{du}{(u-x)v} }.
\end{equation*}
Note again that the contour of integration is a linear combination of the generators ${\mathcal A}, \,{\mathcal B}$ of the homology group of the torus $\gamma=c_1\mathcal A+c_2\mathcal B$ and thus we again obtain a family of solutions to PVI$(\frac{1}{8}, -\frac{1}{8}, \frac{1}{8}, \frac{3}{8})$ parametrized by $c_1/c_2$.

Derivatives of this solution can again be easily computed using the equation of the Schlesinger system which, in this case, are:
\begin{equation*}
\frac{da_1}{dx} =  \frac{a_2}{2x}, \qquad \frac{da_2}{dx} = \frac{a_3-a_2}{2(x-1)}, \qquad a_1+a_2+a_3=0.
\end{equation*}
In particular, with $a_i$ given by \eqref{Picard-Schlesinger}, we have
\begin{equation*}
y'(x) = \frac{xa_1^2 + (x-1)a_2^2 + 2(x-1)a_1a_2}{2(xa_1+(x-1)a_2)^2} \;\;\;\;\;\mbox{ and }\;\;\;\;\;
y''(x) = \frac{\frac{a_2^3}{x} + xa_3^3-2a_1a_2a_3}{2(xa_1+(x-1)a_2)^3}\;.
\end{equation*}

\vskip 0.5cm

We note that our solution to the Schlesinger system from Section \ref{sect_Schlesinger2} in the genus one case coincides with the known solution, see \cite{book}, p. 148. We quote here the derivation from \cite{book} for completeness.

In \cite{book}, the general triangular reduction of the Schlesinger system corresponding to the elliptic curve $v^2=u(u-1)(u-x)$ is considered. Namely, solutions to the system are represented in the form
\begin{equation*}
A^{(1)} = \left(\begin{array}{cc}\alpha & p \\0 & -\alpha\end{array}\right),\qquad
A^{(2)} = \left(\begin{array}{cc}\beta & q \\0 & -\beta\end{array}\right), \qquad
A^{(3)} = \left(\begin{array}{cc}\gamma & -p-q \\0 & -\gamma\end{array}\right), \qquad
A^{(\infty)} = \left(\begin{array}{cc}-\delta & 0 \\0 & \delta\end{array}\right).
\end{equation*}
with $\delta=\alpha+\beta+\gamma$. Then the Schlesinger system reduces to the following differential equations for the functions $p(x)$ and $q(x)$:
\begin{eqnarray}
\label{systempq}
\left\{\begin{array}{c}p_x=\frac{2}{x}\left( \gamma p + \alpha(p+q) \right) \\
\\
q_x=\frac{2}{x-1}\left( \gamma q +\beta(p+q)  \right)
\end{array}\right..
\end{eqnarray}
This system implies a second order ODE for $p(x)$:
\begin{equation*}
x(x-1)p_{xx} + p_x[2\gamma+2\alpha-1 - (4\gamma+2\alpha-1)x] +4p(\gamma^2+\gamma(\alpha+\beta))=0.
\end{equation*}
This is nothing but the hypergeometric equation
\begin{equation}
\label{hypergeometric_eq}
x(1-x)p_{xx} + p_x[c-(a+b+1)x] -abp=0
\end{equation}
with
\begin{equation*}
a=-2(\alpha+\beta+\gamma)=2\delta, \qquad
b=-2\gamma, \qquad
c=1-2(\alpha+\gamma).
\end{equation*}

Now note that with the same choice as in Section \ref{sect_Schlesinger2},
\begin{equation}
\label{parameters}
\alpha=\frac{1}{4}, \qquad
\beta=-\frac{1}{4}, \qquad
\gamma=-\frac{1}{4}, \qquad
\delta=\frac{1}{4},
\end{equation}
we obtain for $p(x)$ the hypergeometric equation of the form $x(x-1)p_{xx}+(2x-1)p_x+\frac{1}{4}p=0$.  It is known (see for example \cite{Clemens}, formula (2.25), p. 61) that this equation is the Picard-Fuchs equation of the corresponding family of elliptic curves and its solutions are linear combinations of periods of the differential $du/v$ on the curves. Thus the corresponding Schlesinger system has solutions of the form
\begin{equation*}
p(x)=c_1\oint_{\mathcal A}\frac{du}{v}+c_2\oint_{\mathcal B}\frac{du}{v},
\end{equation*}
\begin{equation}
\label{q}
q(x)={xc_1}\oint_{\mathcal A}\frac{du}{(u-x)v}+x{c_2}\oint_{\mathcal B}\frac{du}{(u-x)v},
\end{equation}
where $q$ is expressed in terms of $p$ and $p_x$ using the first equation of system \eqref{systempq}. In the case of parameters \eqref{parameters}, this relation is simply $q=2xp_x.$
Note that this solution coincides with \eqref{Picard-Schlesinger}, namely $p(x)=a_1(x), \; q(x)=a_2(x)$ and $-p(x)-q(x)=a_3(x)$.

\begin{remark}
Analogously to the Picard-Fuchs equation for $p(x)$, one obtains that the function $q(x)$ also satisfies a hypergeometric equation as a consequence of system \eqref{systempq}:
\begin{equation*}
x(x-1)q_{xx} + q_x[2(\gamma+\alpha) - (4\gamma+2\alpha+2\beta-1)x] +4q(\gamma^2+\gamma(\alpha+\beta))=0.
\end{equation*}
 The correspondence of parameters with the standard form of the hypergeometric equation \eqref{hypergeometric_eq} is:
$a=2\delta, b=-2\gamma, c=-2(\alpha+\gamma).$
This means that the function $q(x)$ given by \eqref{q} satisfies the following hypergeometric equation:
\begin{equation*}
x(x-1)q_{xx} + 2xq_x +\frac{1}{4}q=0
\end{equation*}
and thus this is one more hypergeometric equation whose solutions are given by periods of Abelian differentials on elliptic curves.
\end{remark}

\section{Known forms of general solution to the Painlev\'e VI  $\left(\frac{1}{8}, -\frac{1}{8}, \frac{1}{8}, \frac{3}{8}\right)$}
\label{sect_others}
For completeness and the convenience of the reader we list known forms of the general solution to the Painlev\'e equation considered here.
For the same family of elliptic curves as above, $v^2=u(u-1)(u-x),$ define $\mu$ to be the period of the curve such that $x=\theta^4_4(0)/\theta^4_4(0)$. Here $\theta_{\bf p,q}(z)=\theta_{\bf p,q}(z|\mu)$ and the well known choices of characteristics  give the Jacobi thetas $\theta_1=-\theta_{\bf 1/2,1/2}$, $ \theta_2=\theta_{\bf 1/2, 0}$, $\theta_3= \theta_{\bf 0,0}$ and $\theta_4=\theta_{\bf 0,1/2}$.

\begin{itemize}

\item In 1987, in \cite{Okamoto} Okamoto obtained the general solution of Painlev\'e VI with parameters (\ref{constants})  by one of his transformations from the following solution of Painlev\'e VI  $\left(0,0,0, \frac{1}{2}\right)$
\begin{equation*}
y_0 (x) = \wp(2c_1 w_1(x) + 2c_2 w_2(x))
\end{equation*}
known as the Picard solution. Here $\wp$ is the Weierstrass function appropriately rescaled such that $u=\wp(z)$ and $v=\wp'(z)$ satisfy the equation of the curve, and $2w_1=\oint_{\mathcal A} \frac{du}{v}$,  $2w_2=\oint_{\mathcal B} \frac{du}{v}$ so that the period of the torus is $\mu=w_2/w_1$.

Given the Picard solution $y_0(x)$,  the general solution $y(x)$ of the
Painlev\'e VI $\left(\frac{1}{8}, -\frac{1}{8}, \frac{1}{8}, \frac{3}{8}\right)$ is obtained by the following Okamoto transformation
\begin{equation}
\label{Okamoto}
y(x) = y_0 + \frac{y_0(y_0-1)(y_0-x)}{x(x-1)y_0^\prime - y_0(y_0-1)}.
\end{equation}
Note that this transformation is written in \cite{Okamoto} in Example 2.1 with a misprint. This is explained and corrected in \cite{DS1}.

\item In 1995, in \cite{Hitchin} Hitchin obtained the following solution:
\begin{equation*}
y(x) = \frac{\theta_1'''(0)}{3\pi^2\theta_4^4(0)\theta'_1(0)} + \frac{1}{3}\left( 1+\frac{\theta_3^4(0)}{\theta_4^4(0)} \right)
%\\
 + \frac{  \theta'''_1(\nu) \theta_1(\nu) - 2\theta''_1(\nu)\theta'_1(\nu) + 4\pi {\rm i} c_1 \left( \theta''_1(\nu)\theta(\nu) - (\theta'_1)^2(\nu) \right)   }{ 2 \pi^2\theta_4^4(0) \theta_1(\nu) \left( \theta'_1(\nu) + 2\pi{\rm i} c_1\theta_1(\nu) \right)}
\end{equation*}
where $\nu=c_1\mu + c_2$.

\item In 1998, in \cite{KiKo} Kitaev and Korotkin obtained the tau-function of a Schlesinger system corresponding to the Painlev\'e equation in question
\begin{equation*}
\tau(x) = \frac{\theta_{\bf p,q}(0)}{\sqrt[8]{x(x-1)}} \left( \int_0^1\frac{du}{v}  \right)^{-\frac{1}{2}}
\end{equation*}
and expressed the solution to the Painlev\'e $\left(\frac{1}{8}, -\frac{1}{8}, \frac{1}{8}, \frac{3}{8}\right)$ depending on two parameters ${\bf p, q}$ in terms of this tau-function as follows:
\begin{equation*}
y(x) = x-x(x-1) \left[ D\left( \frac{\frac{d}{dx}D(\tau)}{\frac{d}{dx}D( \sqrt[8]{x(x-1)}\tau  )}\right)  + \frac{x(x-1)}{D^2\left( \sqrt[8]{x(x-1)}\tau \right)}  \right]^{-1}
\end{equation*}
where $D$ is an operator defined by $D(\cdot) = x(x-1)\frac{d}{dx} {\rm ln}(\cdot).$

\end{itemize}

{\bf Acknowledgements.}
We thank Oleg Lisovyy, Maxim Pavlov and Renat Gontsov for very useful remarks which significantly improved the paper.
The research has been partially supported by the NSF grant 1444147. The research of the first author has been partially
supported by the grant 174020 ``Geometry and topology of manifolds,
classical mechanics, and integrable dynamical systems" of the
Ministry of Education and Sciences of Serbia and by the University
of Texas at Dallas. The second author thanks the Max Planck Institute for Mathematics in Bonn where her part of this work was done for hospitality and support. She also gratefully acknowledges
partial support from the Natural Sciences and Engineering Research Council of Canada and the University of Sherbrooke.

\end{document}